\documentclass[11pt, a4paper]{article}

\usepackage[utf8]{inputenc}
\usepackage[english]{babel}
\usepackage[T1]{fontenc}

\usepackage{amsmath}
\usepackage{amssymb}
\usepackage{amsthm}
\usepackage{mathabx}
\usepackage{tikz}
\usepackage{tikz-cd}
\usepackage{todonotes}
\usepackage{hyperref}
\usepackage{cleveref}
\usepackage{footmisc}
\usepackage{bbm}
\usepackage{afterpage}
\usepackage{mathrsfs}
\usepackage{sseq}
\usepackage{spectralsequences}
\usepackage{extarrows}
\usepackage{enumitem}
\usepackage{stmaryrd}
\usepackage{pdflscape}

\newcommand{\Ell}{\mathrm{Ell}}

\newcommand{\ori}{\mathrm{or}}

\DeclareMathOperator{\Spec}{Spec}
\DeclareMathOperator{\MF}{MF}

\newcommand{\1}{\mathbf{1}}
\newcommand{\BP}{\mathrm{BP}}

\DeclareMathOperator{\MU}{MU}

\newcommand{\Sph}{\mathbf{S}}

\DeclareMathOperator{\TMF}{TMF}
\DeclareMathOperator{\SMF}{SMF}
\DeclareMathOperator{\tmf}{tmf}
\DeclareMathOperator{\Tmf}{Tmf}

\DeclareMathOperator{\Fun}{Fun}

\DeclareMathOperator{\Mod}{Mod}

\DeclareMathOperator{\QCoh}{QCoh}

\DeclareMathOperator{\Sp}{Sp}

\DeclareMathOperator{\Syn}{Syn}

\newcommand{\colim}{\mathrm{colim}\,}

\newcommand{\id}{\mathrm{id}}

\newcommand{\op}{\mathrm{op}}

\newcommand{\E}{\mathbf{E}}

\newcommand{\F}{\mathbf{F}}

\newcommand{\J}{\mathrm{J}}
\renewcommand{\j}{\mathrm{j}}
\newcommand{\M}{\mathsf{M}}

\newcommand{\Z}{\mathbf{Z}}

\newcommand{\al}{\alpha}
\newcommand{\be}{\beta}

\newcommand{\ueq}{\overset{\circ}{=}}

\usepackage{mathtools}

\theoremstyle{theorem}\numberwithin{equation}{section}

\crefname{theorem}{{th}.\!\!}{{ths}.\!\!}
\Crefname{theorem}{{Th}.\!\!}{{Ths}.\!\!}
\newtheorem{theoremalph}{Theorem}

\crefname{theoremalph}{{th}.\!\!}{{ths}.\!\!}
\Crefname{theoremalph}{{Th}.\!\!}{{Ths}.\!\!}
\newtheorem{coralph}[theoremalph]{Corollary}

\crefname{coralph}{{cor}.\!\!}{{cors}.\!\!}
\Crefname{coralph}{{Cor}.\!\!}{{Cors}.\!\!}

\Crefname{problem}{{Prb}.\!\!}{{Prbs}.\!\!}
\newtheorem{prop}[equation]{Proposition}
\Crefname{prop}{{Pr}.\!\!}{{Prs}.\!\!}
\newtheorem{lemma}[equation]{Lemma}
\Crefname{lemma}{{Lm}.\!\!}{{Lms}.\!\!}
\newtheorem{cor}[equation]{Corollary}
\Crefname{cor}{{Cor}.\!\!}{{Cors}.\!\!}

\Crefname{conjecture}{{Conj}.\!\!}{{Conjs}.\!\!}

\theoremstyle{definition}\numberwithin{equation}{section}
\newtheorem{mydef}[equation]{Definition}
\Crefname{mydef}{{Df}.\!\!}{{Dfs}.\!\!}

\Crefname{variant}{{Var}.\!\!}{{Vars}.\!\!}

\Crefname{recall}{{Rcl}.\!\!}{{Rcls}.\!\!}

\Crefname{construction}{{Con}.\!\!}{{Cons}.\!\!}

\Crefname{ass}{{As}.\!\!}{{As}.\!\!}
\newtheorem{question}[equation]{Question}
\Crefname{question}{{Q}.\!\!}{{Qs}.\!\!}

\Crefname{notation}{{Nt}.\!\!}{{Nts}.\!\!}

\Crefname{situation}{{St}.\!\!}{{Sts}.\!\!}

\theoremstyle{remark}\numberwithin{equation}{section}

\Crefname{example}{{Ex}.\!\!}{{Exs}.\!\!}

\Crefname{nonexample}{{NonEx}.\!\!}{{NonEx}.\!\!}

\Crefname{claim}{{Clm}.\!\!}{{Clms}.\!\!}
\newtheorem{remark}[equation]{Remark}
\Crefname{remark}{{Rmk}.\!\!}{{Rmks}.\!\!}

\Crefname{idea}{{Id}.\!\!}{{Ids}.\!\!}

\Crefname{warn}{{Warn}.\!\!}{{Warns}.\!\!}

\Crefname{figure}{{Fig.}\!\!}{{Figs.}\!\!}
\Crefname{footnote}{{Fn.}\!\!}{{Fn.}\!\!}

\Crefname{part}{{\textsection}\!\!}{{\textsection}\!\!}
\Crefname{chapter}{{\textsection}\!\!}{{\textsection}\!\!}
\Crefname{section}{{\textsection}\!\!}{{\textsection}\!\!}
\Crefname{subsection}{{\textsection}\!\!}{{\textsection}\!\!}
\Crefname{appendix}{{\textsection}\!\!}{{\textsection}\!\!}

\begin{document}
\title{
Revisiting the $\be_1$-action on the\\ 
    $3$-primary stable homotopy groups of spheres
}
\author{Jack Morgan Davies\footnote{\href{mailto:jack.davies@ruhr-uni-bochum.de}{\texttt{jack.davies@ruhr-uni-bochum.de}}}}
\maketitle

\begin{abstract}
    Let $\be_1$ be the first $3$-torsion class in the stable homotopy groups of spheres in even degree, well-defined up to a sign. Toda showed that $\be_1^5 \neq 0$, whilst $\be_1^6 = 0$. Shimomura generalised this to the $144$-periodic family generated by $\be_1$, written as $\{\be_{1+9s}\}_{s\geq 0}$, by showing that any $5$-fold product $\prod_5 \be_{1+9s} \neq 0$, whilst all $6$-fold products $\prod_6 \be_{1+9s} = 0$. In this article, we give a simple proof of these results as well as some generalisations to other (non)vanishing $144$-periodic families. Along the way, we also refine known detection results for these families. Our tools include $\BP$-synthetic spectra and the well-known Adams--Novikov spectral sequence for the spectrum of topological modular forms at the prime $3$ together with its Adams operations.
\end{abstract}

\setcounter{tocdepth}{1}
\tableofcontents

\section{Introduction}
Amongst the cornerstones of algebraic topology are the stable homotopy groups of spheres $\pi_n \Sph = \colim \pi_{n+k} S^k$, where $n$ is an integer. These abelian groups assemble into a graded ring $\pi_\ast \Sph = \bigoplus \pi_n \Sph$ by composing and suspending maps between spheres. These rich algebraic structures are very complicated---we only know the order of the abelian groups $\pi_n \Sph$ for roughly $n\leq 90$, for example; see \cite{greenbook,iwxpisph}. In this article, we are interested in a small slice of the multiplicative structure of $\pi_\ast \Sph$. To this end, we work solely at the prime $3$ and implicitly invert all other primes. This is to study the action of a particular 3-torsion element $\be_1$, a generator of $\pi_{10} \Sph \simeq \F_3$ and the smallest even-dimensional class in $\pi_\ast\Sph$, on $\pi_\ast \Sph$ by multiplication. In short, we show that for a host of known nonzero classes in $\pi_\ast \Sph$, multiplication by $\be_1$ is zero, and for another related collection of nonzero classes, multiplication by $\be_1$ is nonzero. These statements would be interesting even if we knew the precise structure of the groups $\pi_n \Sph$ in those degrees containing these products, but these results are obtained without such information.

The classes we are interested in here arise from \emph{chromatic homotopy theory}; a methodology for filtering the ring $\pi_\ast \Sph$ by \emph{chromatic height}. For example, the only degree $\pi_n \Sph$ with nontrivial height $0$ information is $\pi_0 \Sph \simeq \Z$, a consequence of Hopf's computation $\pi_n S^n \simeq \Z$ in \cite{hopfpin} and Serre's finiteness theorem \cite{serreclasses}. Adams' work on the $J$-homomorphism and topological $K$-theory \cite{adamsjofx} computes all of the height $1$ information in $\pi_n \Sph$ for all $n$; also see Toda's \cite{toda_survey}. Although there has been much work on computations at higher heights, perhaps most notably by Miller--Ravenel--Wilson \cite{MRW77} and Goerss--Henn--Mahowald--Rezk \cite{ktwolocalsphere}, height $2$ remains the edge of our knowledge.

The class $\be_1$ can alternatively be defined as the first class of height $2$ in $\pi_\ast \Sph$. The machine of chromatic homotopy theory leads to many other classes of height $2$, the so-called \emph{divided $\beta$-family}. For us, these are classes $\be_{i/j} \in \pi_{16i-4j-2} \Sph$ for various values of $i$ and $j$, and we write $\be_i = \be_{i/1}$. More specifically, works of Behrens--Pemmaraju \cite{marksatya} and Belmont--Shimomura \cite{tmfthree} show the existence of the families of $3$-torsion elements in \Cref{tab:elements} that we focus on here.
\begin{table}[h]
\centering
\begin{tabular}{|l|c|c|c|c|c|c|c|}\hline
Degree mod 144  & 10            & 26            & 37                    & 74            & 81 & 82                & 109                   \\ \hline
Families        & $\be_{1+9s}$  & $\be_{2+9s}$  & $[\al_1\be_{3+9s/3}]$ & $\be_{5+9s}$  & $x_{81}^{(s)}$&$\be_{6+9s/3}$    & $[\al_1\be_{7+9s}]$   \\ \hline
\end{tabular}
\caption{Main characters in the divided $\be$-family used here; square brackets $[xy]$ indicate classes detected by a product $xy$ on the $E_2$-page of the Adams--Novikov spectral sequence, but which no longer a product in $\pi_\ast \Sph$. The family $x_{81}^{(s)}$ for $s\geq 0$ begins with the Toda bracket $x_{81}^{(0)} = \langle \al_1, \al_1, \beta_5\rangle$ and is shown to lift to a $144$-periodic family in \cite[Df.3.28]{heighttwojat3}.}
\label{tab:elements}
\end{table}

We can now state our first theorem concerning the action of any element in the family $\{\be_{1+9s}\}_{s\geq 0}$ on various products of other elements from \Cref{tab:elements}.

\begin{theoremalph}\label{main:equivalences}
    Let $s_i,t\geq 0$ be a collection of nonnegative integers. Then
    \begin{enumerate}
        \item $\prod_{i=1}^A \be_{1+9s_i} \neq 0$ if and only if $A \leq 5$,
        \item $\be_{2+9t}\prod_{i=1}^B \be_{1+9s_i} \neq 0$ if and only if $B \leq 2$, 
        \item $\al_1\be_{2+9t}\prod_{i=1}^C \be_{1+9s_i} \neq 0$ if and only if $C \leq 2$, and
        \item $[\al_1\be_{3/3}]\prod_{i=1}^D \be_{1+9s_i} \neq 0$ if and only if $D \leq 1$.
    \end{enumerate}
\end{theoremalph}

This result is very much inspired by Shimomura's work \cite{shimbetaoneactionatthree}, which states versions of parts 1 and 2 above. Parts 2 and 3 are related, in that the vanishing statement of 2 implies that of 3, and the nonvanishing of 3 implies that of 2. Each item in \Cref{main:equivalences} really consists of two statements: one concerning the vanishing of a large product, and the other concerning the nonvanishing of a smaller product. To prove the vanishing half, we periodify classical differentials in the \emph{Adams--Novikov spectral sequence} (ANSS) for $\Sph$ killing products in low degrees, to obtain differentials killing the periodified products in higher degrees. On the other hand, our nonvanishing arguments use Hopkins' spectrum of \emph{topological modular forms} $\TMF$ together with some cohomology operations to detect this nonvanishing. Both of these halves admit further extensions to the (non)vanishing of other families of products, although the results are not as sharp as \Cref{main:equivalences}; for example, we know that $\al_1 \prod_{i=1}^N \be_{1+9s_i}$ is nonzero for $N\leq 2$ and vanishes for $N\geq 4$, but our techniques do not apply to $N=3$, which is why this case does not appear in \Cref{main:equivalences}. For this reason, we record some generalised vanishing and nonvanishing statements separately. The first we state as a corollary, as it follows immediately by adapting the proof of \Cref{main:equivalences} without any extra computations.

\begin{coralph}\label{main:vanishing}
    Let $s_i,t\geq 0$ be a collection of nonnegative integers. Then the following products of elements in $\pi_\ast\Sph$ all vanish:
    \[\al_1 \prod_{i=1}^4 \be_{1+9s_i}, \qquad \be_{5+9t} \prod_{i=1}^4 \be_{1+9s_i}, \qquad \al_1\be_{6/3} \prod_{i=1}^3 \be_{1+9s_i}\]
\end{coralph}

For our generalised nonvanishing statement, recall that $\TMF$ has an endomorphism $\psi^2$, a kind of \emph{Adams operation}, and that we define $\J^2$ as the equaliser of $\psi^2$ and the identity acting on $\TMF$, a kind of ``height $2$ image-of-$J$ spectrum'' à la Mahowald \cite{mahowaldimageofjspectrum}. By definition, we have maps of ring spectra $\Sph \to \J^2 \to \TMF$.

\begin{theoremalph}\label{main:nonvanishing}
    Let $s_i,t,w\geq 0$ and $x_i \in \{\be_{1+9s}, \be_{6+9s/3}\}_{s\geq 0} \subseteq \pi_{\ast} \Sph$. Then the classes
    \begin{equation}\label{eq:line1}\prod_{i=1}^5 x_i, \qquad \be_{2+9t}\prod_{i=1}^2 x_i, \qquad \be_{5+9t}\prod_{i=1}^2 x_i,\end{equation}
    \begin{equation}\label{eq:line2}\al_1\be_{2+9t}\prod_{i=1}^2 \be_{1+9s_i}, \quad \al_1\be_{2+9t}\prod_{i=1}^2 \be_{6+9s_i/3}, 
    \end{equation}
    \begin{equation}\label{eq:line2.5}
    [\al_1\be_{7+9t}]\prod_{i=1}^2 \be_{1+9s_i}, \quad    [\al_1\be_{7+9t}]\prod_{i=1}^2 \be_{6+9s_i/3}, \quad 
    \beta_{1+9s}[\al_1\be_{3+9t/3}] \be_{6+9w/3},
    \end{equation}
    \begin{equation}\label{eq:line3} x_{81}^{(t)} \be_{6+9w/9} \prod_{i=1}^3 \be_{1+9s_i}, \qquad x_{81}^{(t)} \be_{1+9w} \prod_{i=1}^3 \be_{6+9s_i/3}\end{equation}
    which all vanish in $\J^{2}$ and $\TMF$, do not vanish in $\pi_\ast\Sph$.
\end{theoremalph}

This last statement complements our results of \cite{heighttwojat3} together with Carrick, where we computed the image of $\pi_\ast \Sph$ in $\pi_\ast \J^2$ up to a small discrepancy; one of the main ingredients used here is this \emph{Hurewicz image} of $\J^2$. Combining \Cref{main:nonvanishing} together with \cite[Th.B]{heighttwojat3} gives a catalogue of the known nonvanishing results between products in the divided $\beta$-family at the prime $3$ using $\J^2$.

These results are all obtained using \emph{$\BP$-synthetic spectra} à la Pstr\k{a}gowski \cite{syntheticspectra}, an incarnation of an $\infty$-category of Adams--Novikov spectral sequences. For our vanishing statements, this is mostly an aesthetic choice. Similar arguments can be made explicitly using the ANSS, although our use of synthetic spectra avoids the use of any subtle geometric boundary theorems (such as \cite[Th.2.3.4]{greenbook}). For the nonvanishing statements, the flexibility granted by synthetic spectra, especially in giving us a \emph{modified ANSS} for $\J^2$, is incredibly useful; to the best of our knowledge, the actual ANSS for $\J^2$ has yet to be computed. Moreover, this modified ANSS is simple, see \Cref{pr:manssforJ}, given the ANSS for $\TMF$, which we review in \Cref{sec:tmfsection}. The nonvanishing results above are then a corollary of this modified ANSS, the Hurewicz image of $\J^2$ studied in \cite{heighttwojat3}, and some new explicit arguments.

As an aside, one can study the \emph{mod $\tau$-Hurewicz image} of $\TMF$, so roughly speaking, the image of the $E_2$-page of the ANSS for $\Sph$ in the $E_2$-page of the ANSS for $\TMF$. Simple detection arguments for the ANSS of $\TMF$ as well as the modified ANSS for $\J^2$ (\Cref{pr:modtau_hurewicz_image}) show that the latter detects precisely the divided $\beta$-families
    \[\beta_{1+9t}, \quad \beta_{2+9t}, \quad \beta_{3+9t/3}, \quad \beta_{4+9t}, \quad \beta_{5+9t}, \quad \beta_{6+9t/3}, \quad \beta_{7+9t}, \quad \beta_{8+9t}\]
and no more, refining similar detection results of Belmont--Shimomura \cite{tmfthree} and Carrick and the author \cite{heighttwojat3}; see \Cref{cor:classification}. This result is not just aesthetically pleasing, but from it one can also reprove Shimomura's nonexistence \cite[Ths.F-G]{l2atprime3} of the divided $\beta$-familes $\beta_{3+9t/3}$, $\beta_{4+9t}$, $\beta_{7+9t}$, and $\beta_{8+9t}$; see \Cref{cor:nonexistence}. We hope that these results can aid in future ANSS computations for $\Sph$.

Working with only $\TMF$ and $\J^{2}$, as well as their synthetic versions, does have its limitations. For example, $\be_2^2=0$ in $\J^{2}$, see \Cref{rmk:boundaryissquarezero}, so we cannot recover the fact that $\be_2^2\neq0$ in $\Sph$; see \cite[Tab.A3.4]{greenbook}. The same goes for $\be_5^2$ and the divided $\beta$-families generated by these elements. Moreover, products such as $\al_1 \beta_1^3 \beta_2 \beta_{6/3}$ certainly vanish by \Cref{main:equivalences}, but it is not clear if $\al_1 \beta_1^2 \beta_2 \beta_{6/3}$ vanishes in $\Sph$ despite the fact it vanishes in $\J^2$. In terms of $E_2$-pages, we have $\al_1 \beta_8 \neq 0$ in the ANSS for $\Sph$, but this product vanishes in our modified ANSS for $\J^2$. See \Cref{question:products} and \Cref{rmk:morefamilies?} for more discussion of these limitations and unknowns.

The $Q(N)$-spectra of Behrens' are a further refinement of $\TMF$ and $\J^{2}$, living even closer to $\Sph$; see \cite{buildingandec}. We hope to come back to a discussion of the obvious synthetic version of $Q(N)$, as in \cite[Rmk.2.17]{heighttwojat3}, in the future, utilising either $Q(2)$ or $Q(7)$ at the prime $3$; the former has been used to great success by Behrens' \cite{ktwospheremark}, and the latter would build upon the computations of $\pi_\ast \TMF_0(7)$ of Meier--Ozornova \cite{meieroroztmfo7}. 
This story of the $Q(N)$-spectra is also interesting at the prime $2$, as shown in \cite{markkyle,v232families}, and at primes $p\geq 5$, as shown in \cite{dividedbetafamily}. An optimist might also hope to access information in $\pi_\ast \Sph$ at heights $h\geq 3$ using higher real $K$-theories, perhaps using a model in synthetic spectra via \cite[Th.4.3]{osyn} or \cite{syntheticslicey}, or the \emph{topological automorphic forms} spectra of Behrens--Lawson \cite{taf}. A highly coherent moduli interpretation and equivariant approach will also appear in forthcoming work with Balderrama--Linskens \cite{geometricnorms}.

To add to this speculation, we hope that the results of this article could help us understand the homotopy groups of $\1[\be_1^{-1}]$, the synthetic spectrum $\1$ associated with the ANSS for the sphere $\Sph$ with the class $\be_1$-inverted. Of course, Nishida's nilpotence theorem (or \Cref{main:equivalences}) tells us that classically $\Sph[\be_1^{-1}]=0$, but the synthetic spectrum $\1[\be_1^{-1}]$ is not zero, and captures $\be_1$-periodic behaviour in the ANSS for $\Sph$. Approximations to $\1[\be_1^{-1}]$ are also easy to compute: using (\ref{eq:modtauhomotopyoftmf}) and  \Cref{pr:manssforJ} one easily obtains
\[\pi_{\ast,\ast} \TMF_\BP[\beta_1^{-1}] / \tau \simeq \F_3[\al_1, \beta_1^{\pm}, \Delta^{\pm}]] / (\al_1^2), \quad \pi_{\ast,\ast} \J^2_\BP[\beta_1^{-1}] / \tau \simeq \F_3[\al_1, \beta_1^{\pm}, \Delta^{\pm}, \zeta] / (\al_1^2, \zeta^2),\] 
where $\zeta = \partial 1$ has bidegree $(-1,1)$. The analogous computations without killing $\tau$, as well as their connective variants $\tmf_\BP[\beta_1^{-1}]$ and $\j^2_\BP[\beta_1^{-1}]$, are also straightforward, but a little too cumbersome to write here. This should be the first step in a series of \emph{exotic periodicities} at the prime $3$. At the prime $2$, the $\eta$-inverted synthetic sphere $\1[\eta^{-1}]$ has been computed by Andrews--Miller \cite{invertingeta}. Further work of Andrews \cite{andrewsw1}, Gheorghe \cite{exoticktheories}, and recently by Isaksen--Kong--Li--Ruan--Zhu \cite{w1periodicityisaksenetal} interpret $\eta=w_0$ as the first in a sequence of periodicities $w_n$, and explicitly compute with the ``height $1$ exotic periodicity'' $w_1^4$ which is closely related to $\bar{\kappa} \in \pi_{20,4} \TMF_\BP$. Given our current knowledge of this $\be_1$-action on $\1$, a computation of $\pi_{\ast,\ast}\1[\be_1^{-1}]$ may be within reach. 

\subsection*{Outline}
In \Cref{sec:moorespectra}, we prove the vanishing half of \Cref{main:equivalences} and \Cref{main:vanishing}. To do this, we take differentials in the ANSS and produce $v_2^9$-periodic versions using the language of self-maps of synthetic Moore spectra. In \Cref{sec:tmfsection}, we review some basic facts of $3$-local topological modular forms and $\J^{2}$, including a description of a modified ANSS for $\J^2$, and prove some new facts about its mod $\tau$-Hurewicz image. In \Cref{sec:detection}, we prove the nonvanishing half of \Cref{main:equivalences} and \Cref{main:nonvanishing} using this modified ANSS.

\subsection*{Notation}
As a reminder, everything is implicitly localised at the prime $3$.

We use the notation and foundational facts concerning \emph{$\BP$-synthetic spectra} $\Syn$ of \cite{syntheticspectra}; everything we will use is also laid out in \cite{syntheticsven} in the language of filtered spectra. In particular, $\Syn$ is a stable symmetric monoidal $\infty$-category equipped with a lax monoidal functor $\nu\colon \Sp \to \Syn$ called the \emph{synthetic analogue}, which is a section to the strong monoidal localisation $\tau^{-1}\colon \Syn \to \Sp$ called \emph{$\tau$-inversion} (\cite[\textsection4]{syntheticspectra}). We write $\sigma\colon \Syn \to \Fun(\Z^\op, \Sp)$ for the lax monoidal functor to \emph{filtered spectra} called the \emph{signature} (\cite[Not.2.5]{smfcomputation}), and implicitly use the fact that the signature of $\nu X$ is the ($3$-local) ANSS for $X$ (\cite[Pr.1.25]{osyn}). Following \cite[\textsection2.2]{syntheticj}, we will call $\sigma(X)$ the spectral sequence associated with $X$ and also the (chosen) modified Adams--Novikov spectral sequence for $\tau^{-1}X$. Moreover, we follow the \emph{stem--filtration} grading for $\Syn$, so $\pi_{s,f}$ corresponds to $(s,f)$ in an AN-chart. Formally, we write $\Sigma^{s,f} \1 = \Sigma^{-f}\nu\Sph^{s+f}$ where $\1$ is the unit of $\Syn$. In particular, the element $\tau$ lives in $\pi_{0,-1} \1$ and the $\infty$-categorical suspension has bidegree $(1,-1)$.

We implicitly use the fact that the ANSS for $\1$ and $\TMF$ are both concentrated in those bidegrees $(s,f)$ with $4|s+f$; this is the usual sparsity in the ANSS at $p=3$. We also use explicit computations in these ANSSs, for which we refer the reader to \cite{greenbook} (also see the chart made available by Belmont \cite{evachart}) and \cite[Fig.A.1]{smfcomputation}, respectively. We also fix synthetic lifts of the divided $\beta$-family elements of \Cref{tab:elements}, except for the family $\be_{1+9s}$ for which we fix a specified synthetic lift in \Cref{sec:moorespectra}. In particular, these synthetic lifts have mod $\tau$-reduction given by the Greek letter construction of \cite{MRW77} and \cite[Def.1.3.19]{greenbook}.

\subsection*{Acknowledgements}
Thank you to Adela Zhang and Gabriel Angelini-Knoll for their invitations to Copenhagen and Paris, respectively, and for providing the opportunity for these ideas to come together. Thank you also to Christian Carrick, for our sustained and enjoyable collaboration which led me in this direction and comments on a draft, and to Gabriel, William Balderrama, Eva Belmont, Lennart Meier, Jackson Morris, Sven van Nigtevecht, Liz Tatum, and Jana Wehlburg for some helpful conversations and advice. I would like to heartily thank the anonymous referees for their suggestions and corrections, both of which helped in improving this article. This article was partially written while I was at the Max-Planck Institute for Mathematics in Bonn and was conducted in the framework of the DFG-funded research training group GRK 2240: Algebro-Geometric Methods in Algebra, Arithmetic and Topology.

\section{Vanishing of products in periodic families}\label{sec:moorespectra}
The goal of this section is to obtain the upper bounds from \Cref{main:equivalences}, so to obtain its vanishing statements, and to prove \Cref{main:vanishing}. To do this, we show that chosen synthetic lifts of these elements to $\1$ are $\tau$-power torsion, ie, are hit by differentials in the ANSS for $\Sph$, and hence vanish in $\pi_\ast \1[\tau^{-1}] \simeq \pi_\ast \Sph$. We begin by defining the synthetic divided $\beta$-family $\be_{1+9s}$; this is merely bookkeeping. Then we compute some bigraded homotopy groups of the synthetic Moore spectrum $\1/(3,v_1)$, which then lead us to the desired upper bound.

Start with the mod $3$ Moore spectrum $\Sph/3$. This has a $v_1$-self-map $v_1\colon \Sph^4/3 \to \Sph/3$ by \cite{adamsjofx}; also see \cite[Th.6.3]{syntheticj}. Behrens--Pemmaraju \cite{marksatya} show that the cofibre of this map $\Sph/(3,v_1)$ has a $v_2^9$-self-map of degree $144$. These self-maps yield self-maps of synthetic Moore spectra by applying the synthetic analogue functor
\[v_1\colon \1^{4,0}/3 \to \1/3, \qquad\qquad v_2^9 \colon \1^{144,0}/(3,v_1) \to \1/(3,v_1);\]
one identifies $\1/3 \simeq \nu(\Sph/3)$ and $\1/(3,v_1) \simeq \nu(\Sph/(3,v_1))$ courtesy of \cite[Lm.4.23]{syntheticspectra}. Let us also use the following notation for the maps in the cofibre sequences
\[\1 \xrightarrow{3} \1 \xrightarrow{q_0} \1/3 \xrightarrow{\partial_0} \1^{1,-1}, \qquad \1^{4,0}/3 \xrightarrow{v_1} \1/3 \xrightarrow{q_{1}} \1/(3,v_1) \xrightarrow{\partial_1} \1^{5,-1}/3.\]
Lastly, recall that the right unit formula in the Hopf algebroid $(\BP_\ast, \BP_\ast \BP)$ takes the form
\begin{equation}\label{eq:right_unit}\eta_R(v_2) = v_2 + v_1t^3_1 - v_1^3 t_1,\end{equation}
see \cite[(6.4.16)]{greenbook}, meaning that $v_2$ is primitive modulo $(3,v_1)$, which yields a well-defined element $v_2 \in \pi_{16,0} \1/(\tau,3,v_1)$. It is easy to check that this class admits an abusively named lift $v_2 \in \pi_{16,0} \1/(3,v_1)$.

\begin{mydef}
    For each $s\geq 0$, write $\be_{1+9s} = \partial_0\partial_1(v_2^{9s}\circ v_2) \in \pi_{10+144s,2} \1$.
\end{mydef}

Synthetic lifts for the other families of \Cref{tab:elements} are obtained similarly; one just repeats the constructions of \cite{marksatya,tmfthree} synthetically. This will come up in \Cref{pr:modtau_hurewicz_image,lm:filtrationtwo_stepping_stone}, for example.

Now we can start with our vanishing results. These kinds of arguments can also be made directly at the level of spectral sequences, but since we will need the synthetic language for \Cref{sec:tmfsection,sec:detection}, we have chosen to introduce and use it already.

\begin{lemma}\label{lm:60}
	The homotopy group $\pi_{66,2} \1/(3,v_1)$ vanishes, so in particular, $\tau^8 \be_1^5 v_2 = 0$.
\end{lemma}

\begin{proof}
    The desired homotopy group sits in a short exact sequence
    \[0 \to (\pi_{66,2} \1/3)/v_1 \xrightarrow{q_1} \pi_{66,2} \1 / (3,v_1) \xrightarrow{\partial_1} (\pi_{61,3} \1/3)[v_1] \to 0;\]
    the notation on the left-hand side indicates the $v_1$-quotient on homotopy groups, and the notation on the right-hand side is the $v_1$-torsion on homotopy groups. The right-hand term vanishes as $\pi_{61,3}\1$ and $\pi_{60,4}\1$ both vanish. To see that the left-hand term $(\pi_{66,2} \1/3)/v_1$ also vanishes, start with the computations $\pi_{66,2}\1=0$ and $\pi_{66,2}\1/3 \simeq \F_3$, the latter generated by a lift $x$ of $\tau^4\al_1\be_1\be_2^2$ through the boundary map in the exact sequence
    \begin{equation}\label{eq:boundary_in_proof}0=\pi_{66,2}\1 \to \pi_{66,2} \1/3 \xrightarrow{\partial_0} \pi_{65,3} \1 \to 0.\end{equation}
	However, by \cite[Lm.2.6]{heighttwojat3}, we have $\partial_0(v_1 \cdot q_0(f)) = \al_1 f$ for any element $f\in \pi_{\ast,\ast} X$ of a (3-local) synthetic spectrum $X$. In particular, we see that
	\[\partial_0( v_1 \cdot q_0(\tau^4\al_1\be_1\be_2^2)) = \tau^4\al_1\be_1\be_2^2 = \partial_0(x),\]
	which as $\partial_0$ is injective in this degree, see (\ref{eq:boundary_in_proof}), shows that this lift $x$ is divisible by $v_1$, hence it vanishes in $(\pi_{66,2} \1/3)/v_1$ and this group consequently dies.
\end{proof}

\begin{lemma}\label{lm:56}
    The class $\tau^4 \be_1^2 \be_2 v_2 \in \pi_{62,2} \1/(3,v_1)$ vanishes.
\end{lemma}

\begin{proof}
	The cofibre sequence defining $\1/(3,v_1)$ induces a commutative diagram of abelian groups
    \begin{equation}\label{eq:5lemma}\begin{tikzcd}
        {0}\ar[r]	&	{(\pi_{s,f} \1/3)/v_1}\ar[r, "{q_1}"]\ar[d, "{\cdot \be_1}"] &   {\pi_{s,f} \1/(3,v_1)}\ar[r, "{\partial_1}"]\ar[d, "{\cdot \beta_1}"]  &   {(\pi_{s-5,f+1} \1/3)[v_1]}\ar[d, "{\cdot \be_1}"]\ar[r]	&	{0}\\
        {0}\ar[r]	&	{(\pi_{s+10,f+2} \1/3)/v_1}\ar[r, "{q_1}"] &   {\pi_{s+10,f+2} \1/(3,v_1)}\ar[r, "{\partial_1}"]   &   {(\pi_{s+5,f+3}\1/3)[v_1]}\ar[r]	&	{0}
    \end{tikzcd}\end{equation}
    with short exact rows. As $\tau^4 \beta_1^3 \beta_2 = 0 \in \pi_{56,8} \1$, this class is hit by a $d_5$ in the ANSS for $\Sph$, it suffices to show that multiplication by $\be_1$ is injective on $\pi_{62,2} \1/(3,v_1)$. From (\ref{eq:5lemma}), it suffices to see that the left-hand and right-hand vertical maps are injective for $(s,f) = (62,2)$. This is clear, because $\pi_{57,3} \1/3 =0$ as $\pi_{57,3}\1 = \pi_{56,4}\1 = 0$, where as $(\pi_{62,2} \1/3)/v_1 \simeq \F_3$ is generated by $\tau^4 \beta_1 \beta_2^2$ and $\tau^4 \beta_1^2 \beta_2^2 \neq 0$ in $(\pi_{72,4}\1/3)/v_1$.
\end{proof}

The same diagram chase also shows the next result, setting $(s,f) = (63,1)$ in (\ref{eq:5lemma}).

\begin{lemma}\label{lm:57}
    The class $\tau^4 \be_1 [\al_1\be_{3/3}] v_2 \in \pi_{63,1} \1/(3,v_1)$ vanishes.
\end{lemma}


\begin{remark}\label{rmk:tautologicalvanishing}
    Inside $\pi_{\ast,\ast}\1$ the products $\tau^4\al_1\be_1^3$, $\tau^4 \be_1^3 \be_5$, and $\tau^4\al_1 \be_1^2 \be_{6/3}$ all vanish. We have not been able to show that either $\tau^4 \al_1 \be_1^2 v_2\in \pi_{39,1} \1/(3,v_1)$, $\tau^4 \be_1^2 \be_5 v_2 \in \pi_{110,2} \1/(3,v_1)$, nor $\tau^4\al_1\be_{6/3}\be_1 v_2\in \pi_{111,1}\1/(3,v_1)$ vanish, meaning that we cannot conclude that the associated classical differentials periodify; this is formulated as \Cref{question:products}. In \Cref{rmk:periodifyingsource}, we briefly discuss a similar problem of trying to periodify the \emph{source} of these differentials. Note that the first three products above vanish in $\1$, then the associated products in $\1/(3,v_1)$ clearly vanish after multiplication by $\be_1$; this leads to the vanishing of the families in (\ref{eq:vanishinglinetwo}) below, and seems worth recording, even though it is obvious.
\end{remark}

A simple argument turns the lemmata above into vanishing statements in $\1$ and $\Sph$.

\begin{cor}\label{cor:vanishinghalf}
	For $s_i,t\geq 0$, the products
	\begin{equation}\label{eq:vanishinglineone}\prod_{i=1}^6 \be_{1+9s_i}, \qquad \be_{2+9t}\prod_{i=1}^3 \be_{1+9s_i},
    \qquad [\al_1 \be_{3/3}]\prod_{i=1}^2 \be_{1+9s_i}\end{equation}
    \begin{equation}\label{eq:vanishinglinetwo}\al_1 \prod_{i=1}^4 \be_{1+9s_i}, \qquad \be_{5+9t} \prod_{i=1}^4 \be_{1+9s_i}, \qquad \al_1\be_{6/3} \prod_{i=1}^3 \be_{1+9s_i}\end{equation}
	are all $\tau^4$-torsion, except for the first $6$-fold product, which is $\tau^8$-torsion, in $\pi_{\ast,\ast}\1$. In particular, all of the associated products in $\pi_\ast \Sph$ vanish.
\end{cor}

\begin{proof}
    Let us start with the family of classes $\prod_{i=1}^6 \be_{1+9s_i}$. Consider the synthetic product $\be_1^5 \be_{1+9s} \in \pi_{60+144s,12} \1$ for some $s\geq 0$. This class is defined as $\be_1^5 \partial_0\partial_1(v_2^{9s}\circ v_2)$. Using the $\1$-linearity of $\partial_0\partial_1$ and $v_2^{9s}$, we obtain
    \[\tau^8 \be_1^5 \partial_0\partial_1(v_2^{9s}\circ v_2) = \partial_0\partial_1(v_2^{9s}\circ(\tau^8 \be_1^5 v_2)) = 0,\]
    the last equality courtesy of \Cref{lm:60}. We will see in \Cref{pr:nonvanishingbeta1beta63}, which is proven independently, that these classes $\be_1^5 \be_{1+9s}$ are \textbf{not} $\tau^4$-torsion, so the mod $\tau$-reduction of these classes are hit by $d_9$-differentials in the ANSS for $\Sph$. Using Toda's result that for any $s+t = u+v$ we have
    \begin{equation}\label{eq:toda_equality}uv\be_s \be_t = st \be_u \be_v \in \pi_{\ast,\ast} \1/\tau,\end{equation}
    see \cite[Th.5.3]{todabetaproducts} for the original statement or \cite[Th.5.6.4]{greenbook} for this form, we obtain
    \[\be_1^5 \be_{1+9s} = \prod_{i=1}^6 \be_{1+9s_i} \in \pi_{\ast,\ast} \1/\tau\]
    as long as $s = \sum_{i=1}^6 s_i$. In particular, we see that every product of the form $\prod_{i=1}^6 \be_{1+9s_i}$ has a mod $\tau$-reduction hit by a $d_9$, hence these products are also $\tau^8$-torsion. Inverting $\tau$ then shows that products of the form $\prod_{i=1}^6 \be_{1+9s_i}$ vanish in $\pi_\ast\Sph$. The proofs of the other two families of (\ref{eq:vanishinglineone}) of products are the same, referring to \Cref{lm:56} and \Cref{lm:57} when necessary. For the families of (\ref{eq:vanishinglinetwo}), we refer to \Cref{rmk:tautologicalvanishing},
\end{proof}

\begin{proof}[Proof of \Cref{main:vanishing}]
    This is just the vanishing of (\ref{eq:vanishinglinetwo}) from \Cref{cor:vanishinghalf}.
\end{proof}

By the nonvanishing results to come, the products in (\ref{eq:vanishinglineone}) do not vanish if we divide them by $\be_{1+9s}$, however, what follows does not answer this question for those products in (\ref{eq:vanishinglinetwo}).

\begin{question}\label{question:products}
	Do the products of (\ref{eq:vanishinglinetwo}) vanish in $\pi_\ast\Sph$ when one factor of $\be_{1+9s}$ is removed?
\end{question}


\section{The modified Adams--Novikov spectral sequence for \texorpdfstring{$\J^2$}{J2}}\label{sec:tmfsection}
In this section, we give an exposition of the $\E_\infty$-ring $\J^2$, the synthetic $\E_\infty$-ring $\J^2_\BP$, and its associated modified Adams--Novikov spectral sequence, ending with some detection statements.

We start with Hopkins' spectrum $\TMF$ of \emph{topological modular forms}, defined as the global sections of the spectral moduli stack $\M_\Ell^\ori$ of \emph{oriented elliptic curves}; see \cite[\textsection7]{ec2} for this definition and \cite{uniqueotop} for the comparison to the Goerss--Hopkins--Miller model of \cite{tmfbook}. Hopkins--Mahowald computed the homotopy groups of $\TMF$ using its ANSS, see \cite[Cor.E]{smfcomputation} or \cite[\textsection5.8]{lennartlecturenotes} for details; one can also compute the $E_2$-page using \cite[\textsection5]{bauer} and then invert $\Delta$. We will use $\nu\TMF$ as our synthetic lift of $\TMF$, and denote it by $\TMF_\BP$. By \cite[Th.C]{osyn}, this synthetic spectrum can be identified with \emph{synthetic modular forms} $\SMF$, which naturally captures the \emph{descent spectral sequence} (DSS) for $\TMF$. In fact, \emph{ibid} identifies the DSS for $\TMF$ with its ANSS. We implicitly use this identification below and write the $E_2$-page of the ANSS for $\TMF$ in terms of the cohomology of the classical moduli stack of elliptic curves $\M_\Ell$. Variants of $\TMF$, such as $\Tmf$ and $\tmf$, are discussed in \Cref{rmk:connective_vs_not}.

\begin{remark}
	Let us not work implicitly localise at the prime 3 for this remark, as this is a subtle point concerning localisations of $\TMF$ that is often swept under the rug. One defines $\TMF$ as the global sections $\Gamma(\M_\Ell^\ori)$, so when one wants to compute the DSS with a set of primes $S$ inverted, one wants to compute the $E_2$-page using the cohomology of $\M_\Ell \times \Spec \Z[S^{-1}]$; in this article we would take $S$ to be the set of all primes excluding $3$. However, this descent spectral sequence converges to $\Gamma(\M_\Ell^\ori \times \Spec \Sph[S^{-1}])$, an object which \emph{a priori} is different to $\TMF[S^{-1}]$. It would then be convenient if the natural map of $\E_\infty$-rings
	\[\TMF[S^{-1}] \to \Gamma(\M_\Ell^\ori, \mathcal{O}_{\M_\Ell^\ori}[S^{-1}]) \simeq \Gamma(\M_\Ell^\ori \times \Spec \Sph[S^{-1}])\]
	were an equivalence, which would remove this potential difference. The global sections functor $\Gamma \colon \QCoh(\M_\Ell^\ori) \to \Mod_{\TMF}(\Sp)$ from quasi-coherent sheaves on $\M_\Ell^\ori$ to $\TMF$-modules in spectra can be computed as a limit, and it is not true that limits generally commute with localisations, a kind of colimit. In this case, we are saved by a famous result of Mathew--Meier \cite{akhilandlennart}, who show that the global sections functor above is an equivalence, hence commutes with colimits; a purely spectral algebro-geometric proof of this fact can be found in \cite[Th.A]{reconstruction}. The map of $\E_\infty$-rings above is then an equivalence, meaning that the spectral sequence we want to compute, actually computes what we want it to compute.
\end{remark}

The computation of the ANSS for $\TMF$ can be encoded synthetically as follows: its $E_2$-page is given by the mod $\tau$ homotopy groups of $\TMF_\BP$
\begin{equation}\label{eq:modtauhomotopyoftmf}
\pi_{\ast,\ast} \TMF_\BP/\tau \simeq \frac{\Z[c_4, c_6, \Delta^{\pm}, \al, \be]}{(c_4^3-c_6^2 - 1728\Delta, 3\al, \al^2, \al c_4, \al c_6, 3\be, \be c_4, \be c_6)}
\end{equation}
generated by elements of bidegree
\[|c_i|=(2i,0),\quad |\Delta|=(24,0),\quad |\al|=(3,1),\quad |\be| = (10,2);\]
the elements $c_4$ and $c_6$ correspond to the normalised Eisenstein series and $\Delta$ to the discriminant modular form following \cite[\textsection III.1, p.42]{silvermaneasy}; for this explicit expression, with slightly different normalisation conventions, see \cite[\textsection5.6 \& 5.8]{lennartlecturenotes}. The elements $\al$ and $\be$ also have algebraic descriptions: in the exact sequence on sheaf cohomology
\[H^0(\M_\Ell, \omega^{\otimes 2}) \to H^0(\M_{\Ell,\F_3}, \omega^{\otimes 2}) \to H^1(\M_\Ell, \omega^{\otimes 2}) \simeq \F_3\{\al\},\]
the class $\al$ is the obstruction to lifting the mod $3$ modular form $b_2$ of weight $2$ to an integral modular form, see \cite[Pr.7.2(II)]{delignecohomologyMELL}, and $\be$ is the Massey product $\langle \al, \al, \al\rangle$. There are only $d_5$'s and $d_9$'s in this spectral sequence \cite[\textsection7.1]{smfcomputation}, multiplicatively generated by
\[d_5(\Delta) = \pm \al\be^2,\qquad d_9([\al\Delta]) = \pm \be^5.\]
These can be deduced by knowing that these classes $\al$ and $\beta$ detect $\al_1$ and $\beta_1$ in $\pi_\ast \Sph$, and the known facts that $\al_1 \beta_1^3$ and $\beta_1^6$ both vanish in $\pi_\ast \Sph$ which where known to Toda; see \cite{toda_an_important_relation} or \cite[Th.4.4.2]{greenbook}. This computes the bigraded homotopy groups of $\TMF_\BP$ as
\[\pi_{\ast,\ast} \nu\TMF \simeq \Z\left[\tau, c_4, c_6, [3\Delta], [c_4\Delta], [c_6\Delta], [3\Delta^2], [c_4\Delta^2], [c_6\Delta^2], \Delta^{\pm 3}, \al, \be, [\al\Delta]\right] / I,\]
where $I$ is a large ideal not worth writing down here; it reflects the facts that 
these homotopy groups yield $\pi_\ast\TMF$ upon inverting $\tau$ and (\ref{eq:modtauhomotopyoftmf}) upon killing $\tau$. We suggest that the reader refers to a figure such as \cite[Fig.A1]{smfcomputation} (and invert $\Delta^3$) to best digest this information; that is certainly how we think about it.

Now we move on to Adams operations on $\TMF$.

\begin{mydef}\label{def:operationsonTMF}
    For an integer $3\nmid k$, let $\psi^k\colon \TMF \to \TMF$ be the \emph{$k$th Adams operation}, the maps of $\E_\infty$-rings of \cite[Df.2.1]{heckeontmf}.
\end{mydef}

\begin{remark}\label{rmk:homotopymodtaucomputation}
By \cite[Cor.2.12]{heckeontmf}, we can compute these operations on $\TMF_\BP$ modulo $\tau$ as the corresponding algebraic operations. In particular, if $f\in \pi_{2d,0} \TMF_\BP/\tau \simeq \MF_d$ is a modular form of weight $d$, then we have $\psi^k(f) = k^d f$. One can also compute $\psi^k(\al) = \al$ and $\psi^k(\be) = \be$ as in $\TMF$ these classes lie in the image of the unit map $\Sph \to \TMF$. As these classes generate $\pi_{\ast,\ast}\TMF_\BP/\tau$ as a ring, see (\ref{eq:modtauhomotopyoftmf}), this gives us a computation of $\psi^k$ on the bigraded homotopy groups of $\TMF_\BP/\tau$, as $\psi^k$ is multiplicative. These kinds of computations of Adams operations on elements in filtration $0$ are true for generic Adams operations produced by chromatic homotopy theory; see \cite[Pr.6.18]{luriestheorem}. Computations of Adams operations on $\TMF$ without killing $\tau$, including at the prime 2, follow from the computations on $\pi_\ast \TMF$ from \cite[Th.C]{adamsontmf}. 
\end{remark}

In \cite{heighttwojat3}, together with Carrick, we studied the equaliser of $\psi^2$ and the identity acting on $\TMF$, first appearing in \cite[Df.4.9]{adamsontmf}, as well as the appropriate synthetic equaliser.

\begin{mydef}
    Write $\J^{2}$ for the $\E_\infty$-ring given as equaliser of the two $\E_\infty$-endomorphisms $\psi^2,\id$ of $\TMF$, or equivalently, the $\E_0$-ring given by the fibre of $\psi^2-1$. Write $\J^{2}_\BP$ for the synthetic $\E_\infty$-ring given as the equaliser of the two maps $\nu\psi^2,\id$ on $\TMF_\BP$. In both cases, we write
    \[\partial \colon \Sigma^{-1}\TMF \to \J^{2}, \qquad \partial \colon \Sigma^{-1,1}\TMF_\BP \to \J_\BP^{2}\]
    for the boundary map of their defining fibre sequences.
\end{mydef}

The ``$\J$''-notation hearkens back to Whitehead's \emph{$J$-homomorphism} and its relation to Adams operations on topological $K$-theory discussed in \cite{adamsjofx}. The Adams operation $\psi^2$ is used for simplicity here; most of this article also works for any $\psi^k$ with $3\nmid k$.

\begin{remark}\label{rmk:connective_vs_not}
    Many of the arguments made here also apply to either \emph{projective topological modular forms} $\Tmf$ or \emph{connective modular forms} $\tmf$. The former has a spectral algebro-geometric interpretation as the global sections of the \emph{compactification of the moduli stack of oriented elliptic curves} $\overline{\M}_\Ell^\ori$, sketched in \cite[\textsection4.3]{lurieecsurvey} and carried out in \cite[Cor.B]{globaltate}, and one defines $\tmf = \tau_{\geq 0} \Tmf$. The DSS for $\Tmf$ and the ANSS for $\tmf$ were computed by Hopkins--Mahowald in the 1990's, but full details of this computation only appeared in \cite{smfcomputation}; see \textsection1.1 of \emph{ibid} for more history and background, and why one never wants to compute the ANSS for $\Tmf$, as well as \cite{osyn} for a synthetic implementation of these spectral sequences Moreover, both $\Tmf$ and $\tmf$ also have Adams operations by \cite[Th.A]{adamsontmf}, although they are not known to be as well-behaved as those on $\TMF$; compare \cite[Ths.C \& D]{heckeontmf} with \cite[Th.C]{realspectra}, for instance. From a computational perspective, Adams operations on $\tmf$ can be easier to work with as the homotopy groups of $\tmf$ are much smaller, hence why its use in \cite{heighttwojat3}. In this article, the only reasons to focus on $\TMF$ are to add computations of $\J^2$ explicitly to the literature and use fact that a power of $\Delta$ acts invertibly on $\TMF$ to more easily prove \Cref{lm:filtrationtwo_stepping_stone} below. 
\end{remark}

The following is a nonconnective version of \cite[Pr.2.12]{heighttwojat3}.

\begin{prop}
    The synthetic spectrum $\J^{2}_\BP$ is a $\tau$-complete and synthetic lift of $\J^{2}$.
\end{prop}

\begin{proof}
    The Hopkin--Mahowald equivalence $\TMF \otimes X_3 \simeq \TMF_0(2)$ for a particular finite spectrum $X_3$, see \cite[(4)]{tmfhomology}, shows that $\TMF$ is $\MU$-nilpotent complete. By \cite[Pr.A.13]{burkhahnseng}, $\TMF_\BP$ is $\tau$-complete, and hence so is $\J^2_\BP$ as $\tau$-complete synthetic spectra are closed under limits. The fact that $\tau$-inversion is exact and is a retract of $\nu$ shows that these synthetic spectra are lifts of the indicated spectra.
\end{proof}

As $\J^{2}_\BP$ is a synthetic lift of $\J^{2}$, its associated signature spectral sequence converges to the homotopy groups of $\J^{2}$. Moreover, the $E_2$-page of this spectral sequence is readily computable from the $E_2$-page of the ANSS for $\TMF$. For degree reasons, there are also not many differentials or extensions; this was used crucially in \cite{heighttwojat3}.

\begin{figure}[h]\begin{center}
\makebox[\textwidth]{\includegraphics[trim={3.5cm 15.5cm 1.5cm 4.2cm},clip,page = 1, scale = 1]{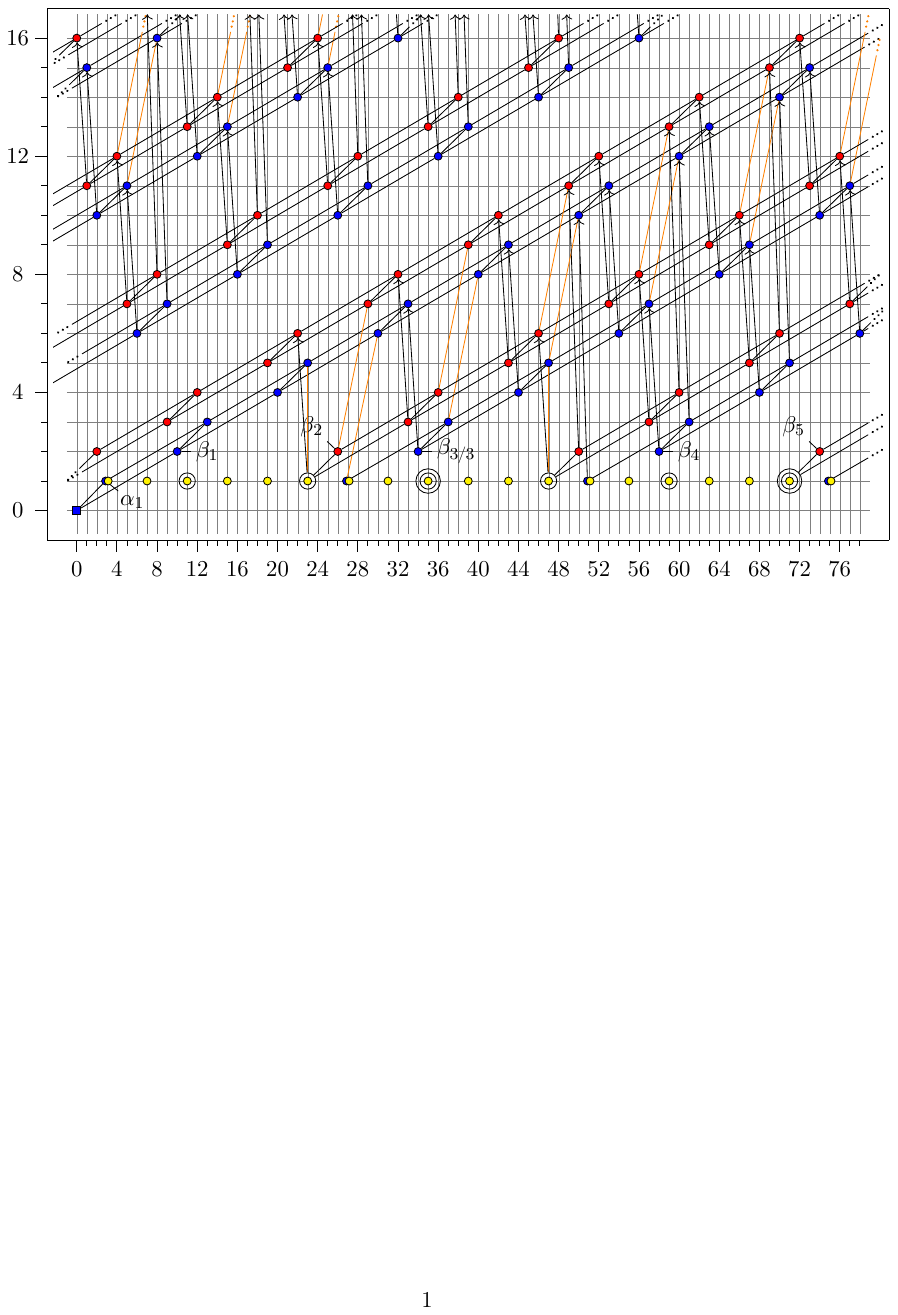}}
\caption{\label{ss1} Signature spectral sequence of $\J^{2}_\BP$ for stems $s$ in $0 \leq s \leq 78$; see \Cref{pr:manssforJ}.}
\end{center}\end{figure}

\begin{figure}[h]\begin{center}
\makebox[\textwidth]{\includegraphics[trim={3.5cm 15.5cm 1.5cm 4.2cm},clip,page = 2, scale = 1]{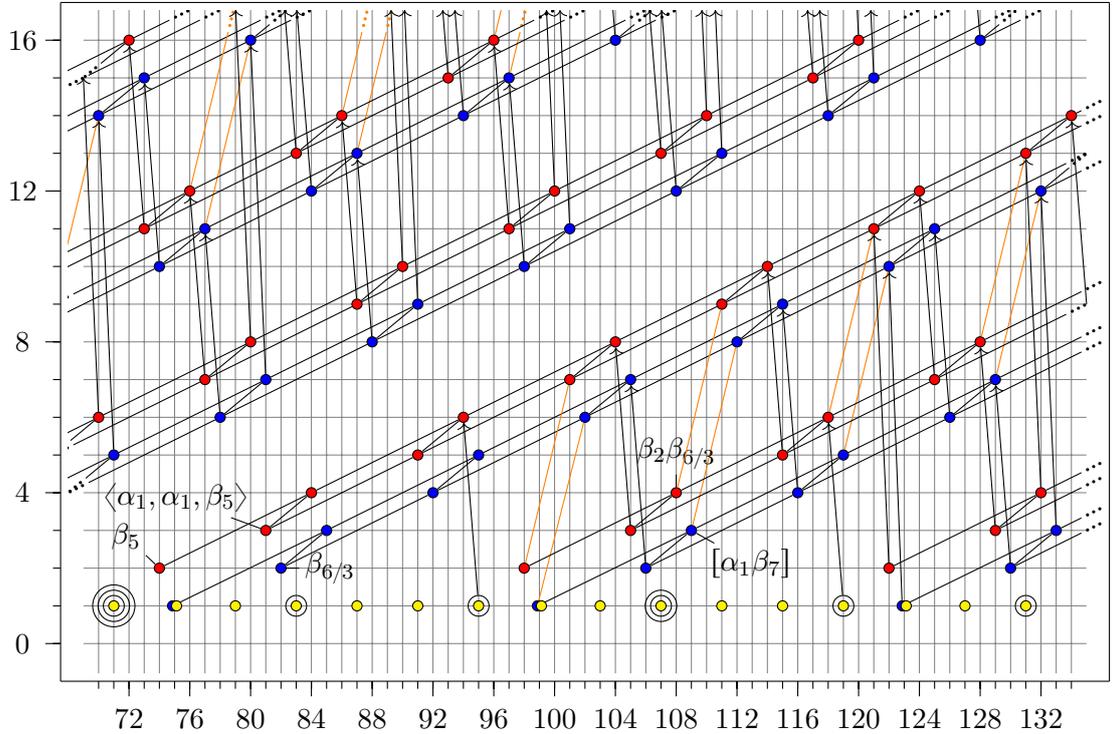}}
\caption{\label{ss2} Signature spectral sequence of $\J^{2}_\BP$ for stems $70\leq 150$; see \Cref{pr:manssforJ}.}
\end{center}\end{figure}

\begin{prop}\label{pr:manssforJ}
    The signature spectral sequence associated with $\J^{2}_\BP$ has differentials determined by the fibre sequence
    \[\Sigma^{-1,1}\TMF_\BP \xrightarrow{\partial} \J^{2}_\BP \xrightarrow{p} \TMF_\BP;\]
    see \Cref{ss1,ss2}. There is one exotic $3$-extension from $[\partial(\Delta^{1+3s})]$ to $[\al_1\be_1^2\Delta^{3s}]$ for each $s\geq 0$. In these charts, blue signifies classes lifted from the ANSS for $\TMF$, and red and yellow mean the image of the ANSS for $\TMF$ along the boundary map; the yellow classes follow an image-of-$J$ pattern. Red and blue dots are copies of $\F_3$, the blue rectangle is a $\Z[j^3]$ for $j=\frac{c_4^3}{\Delta}$, and yellow circles are isomorphic to $\Z[j]/(2^{r}-1)$, where $r=\tfrac{s+1}{2}$ for a yellow dot in the $s$-stem. Lines indicate multiplication by $\al_1$ or $\be_1$; multiplication lines with yellow source are omitted.
\end{prop}

In other words, once we use $\partial$ and $p$ to induce as many differentials as possible, there is no room for any other differentials. This signature spectral sequence also gives the structure of $\pi_{\ast,\ast} \J^{2}_\BP$ by the usual omnibus theorem; see \cite[\textsection A]{burkhahnseng} or \cite[Th.3.62]{syntheticsven}.

\begin{proof}
    For the $E_2$-page, one uses the computation of $\psi^2$ on $\pi_{\ast,\ast}\TMF_\BP/\tau$ described in \Cref{rmk:homotopymodtaucomputation}. The differentials are clear, as once one projects and lifts as many differentials as are in $\TMF_\BP$, there is no more room for any more differentials in $\J^2_\BP$. The exotic $3$-extension comes from the fact that if we compute $\pi_{23+72s}\J^{2}$, then all of the generators can be taken to be simple $9$-torsion, the image of $\pi_{24+72s}\TMF$. This forces these extensions in the spectral sequence.
\end{proof}

\begin{remark}\label{rmk:boundaryissquarezero}
    The product structure on $\J^{k}$ can be made reasonably explicit. The Bousfield--Kan filtration of $\J^{k}$, placing the source $\TMF$ of $\psi^k$ and $\id$ in filtration $0$ and the target $\TMF$ in filtration $1$, is a multiplicative filtration. In particular, two elements in filtration $1$ multiply to elements of filtration $\geq 2$, which is trivial. Specifically, given $x,y\in \pi_\ast\TMF$, we have $\partial(x)\partial(y)=0$ in $\pi_\ast\J^{2}$. The same goes for $\J^{2}_\BP$ using a Bousfield--Kan spectral sequence in synthetic spectra. Visually, this means products of red classes vanish. Products of blue classes with either red or blue classes can usually be explicitly computed in $\TMF$.
\end{remark}

The final piece of information we need about $\J^2_\BP$ is the image of $\pi_{\ast,\ast}\1$ in its homotopy groups, so the \emph{synthetic Hurewicz image} of $\J^2_\BP$. This follows from \cite{tmfthree,heighttwojat3}.

\begin{prop}\label{pr:synthetihurewicz}
    The chosen synthetic lifts of all of the classes in \Cref{tab:elements} have nonzero image in $\pi_{\ast,\ast}\J^2_\BP$. In more detail, for $t\geq 0$, the classes
    \[\be_{1+9t}, \quad \be_{2+9t}, \quad [\al_1\be_{3+9t/3}], \quad\be_{5+9t}, \quad x_{81}^{(t)}, \quad \be_{6+9t/3}, \quad [\al_1\be_{7+9t}]\]
    have image in $\J^2_\BP$ given by the classes
    \[[\be\Delta^{6t}],  \partial([\al\Delta]\Delta^{6t}),  \be_1([\al\Delta]\Delta^{6t}), \partial([\al\Delta^3]\Delta^{6t}]),  \partial(\be\Delta^{3+6t}),  [\be\Delta^{3+6t}], \be_1([\al\Delta]\Delta^{3+6t}),\]
    respectively, up to a unit. In particular, all of the products of \Cref{main:equivalences}, \Cref{main:vanishing} and \Cref{main:nonvanishing} have nonzero image in $\J_\BP^{2}$.
\end{prop}

Of course, having nonzero image in $\J_\BP^2$ does not mean that an element is nonzero in $\J_\BP$; it may be hit by a differential in the modified ANSS for $\J^2$.

\begin{proof}
	This detection statements holds after inverting $\tau$, see \cite[Th.6.5]{tmfthree} and \cite[Th.A]{heighttwojat3}. The fact that $\pi_{\ast,\ast}\J^2_\BP$ is $\tau$-torsion free in these degrees, see \Cref{pr:manssforJ} and \Cref{ss1,ss2}, shows that this detection also must hold synthetically.
\end{proof}

Although it will not play a role in the rest of this article, it is worth noting that the \emph{mod $\tau$-Hurewicz image} of $\TMF_\BP$ and $\J^2_\BP$, so the image of the unit $\1/\tau$, also has some interesting applications to the \emph{nonexistence} of certain divided $\beta$-families; compare with \cite[Th.F]{l2atprime3}.

\begin{prop}\label{pr:modtau_hurewicz_image}
	The families $\beta_{3+9t/3}$, $\beta_{4+9t}$, and $\beta_{7+9t}$ are all nonzero in $\TMF_\BP/\tau$. The family $\beta_{8+9t}$ is nonzero in $\J^2_\BP/\tau$.
\end{prop}

Detection of the $\beta_{4+9t}$-family in $\J^2_\BP/\tau$ also follows from Toda's equality (\ref{eq:toda_equality}):
\[\beta_{1+9a} \beta_{5+9b} \ueq \beta_{2+9c}\beta_{4+9d} \in \pi_{\ast,2} \1/\tau\]

\begin{proof}
	For the $\beta_{3+9t/3}$- and $\beta_{7+9t}$-families, this follows from \cite[Th.6.5]{tmfthree} as products of these classes are detected. For the $\beta_{4+9t}$-family, we can argue directly as in \cite{heighttwojat3} or \cite{syntheticj}; the argument is the same for all $t$ with only superficial differences, so let $t=0$. The class $\beta_4$ can be defined as $\partial_0(\partial_1(v_2^4))$ in $\pi_{58,2} \TMF_\BP/\tau$. First, note that $v_2^4 \in \pi_{64,0} \TMF_\BP/(\tau,3,v_1)$ is nonzero as its image in $\TMF_0(2)_\BP /(\tau,3,v_1)$ is nonzero; this can be seen using the decomposition $\tmf_0(2) \simeq \BP\langle 2\rangle \oplus \Sigma^8 \BP\langle 2 \rangle$ of Hill \cite[Pr.2.3]{hillsigmathree}, where it is credited to Hopkins--Mahowald and Behrens, together with $\tmf_0(2)[\Delta^{-1}] \simeq \TMF_0(2)$. The map
	\[\Z/3 \{c_4^{3r} c_6 \Delta^{2-r}\}_{r\geq 0} \simeq \pi_{60,0} \TMF_\BP/(\tau,3) \xrightarrow{\cdot v_1} \pi_{64,0} \TMF_\BP/(\tau,3) \simeq \Z/3 \{c_4^{3r+2}\Delta^{2-r}\}_{r\geq 0}\]
	is an isomorphism as $c_{2i}$ is congruent to $v_1^{i}$ modulo $3$, so $\partial_1$ is an injection in this degree. As $\pi_{59,1} \TMF_\BP/\tau$ vanishes, $\partial_0$ is also an injection in this degree. In total, we see that $\beta_4$ is nonzero in $\TMF_\BP/\tau$, as desired; it is detected by $\beta \Delta^2$ up to a unit.

    Finally, for the $\beta_{8+9t}$-case in $\J^2_\BP/\tau$, we use the previous facts that the $\beta_{1+9t}$-, $\beta_{4+9t}$-, and $\beta_{5+9t}$-families are all detected in $\J^2_\BP/\tau$, and Toda's equality (\ref{eq:toda_equality}) to obtain
    \[\beta_{1+9a}\beta_{8+9b} \ueq \beta_{4+9c}\beta_{5+9d} \in \pi_{\ast,2} \J^2_\BP/\tau\]
    for any $a,b,c,d\geq 0$ with $a+b=c+d$. As the right-hand product is nonzero in $\J^2_\BP/\tau$, see \Cref{pr:manssforJ}, the $\beta_{8+9t}$-family nonzero in $\J^2_\BP/\tau$, finishing the proof.
\end{proof}

This recovers some nonexistence statements of Shimomura; see \cite[Ths.F-G]{l2atprime3}.

\begin{cor}\label{cor:nonexistence}
	The divided $\beta$-family elements $\beta_{3+9t/3}$, $\beta_{4+9t}$, $\beta_{7+9t}$, and $\beta_{8+9t}$ do not exist in $\pi_\ast \Sph$ for $t \geq 0$, ie, they do not survive the ANSS for $\Sph$.
\end{cor}

\begin{proof}
	All of these classes lie in the mod $\tau$-Hurewicz image of either $\TMF_\BP$ or $\J^2_\BP$ by \Cref{pr:modtau_hurewicz_image}, and they all support differentials in the associated ANSS for $\TMF$ or the modified ANSS for $\J^2$, respectively, hence these differentials must also be nonzero in that of the sphere.
\end{proof}

\begin{remark}\label{rmk:periodifyingsource}
	It is tempting to claim that because the divided $\beta$-families of \Cref{cor:nonexistence} all support differentials in the sphere, that it follows that the targets of these differentials in the ANSS of $\TMF$ must also vanish in the sphere. However, this is not immediately clear without more work. For instance, we can only conclude that $d_5(\beta_{4+90}) = \pm \beta_1^2 [\al \beta_{3+90/3}] + x$ for some class $x$ in the kernel of the unit map from the $E_2$-page of the ANSS of $\Sph$ to that of $\TMF$.
\end{remark}

The mod $\tau$-Hurewicz image of $\TMF_\BP$ is determined by \Cref{pr:modtau_hurewicz_image} and \Cref{lm:filtrationtwo} below. Keeping in mind this uncertainty and the added uncertainty of \cite[Th.A]{heighttwojat3}, we also understand most of the mod $\tau$-Hurewicz image of $\J^2_\BP$. Combined with \Cref{lm:filtrationtwo_stepping_stone} below, yields a neat classification of divided $\beta$-family elements in terms of $\J^2_\BP$; see \Cref{cor:classification}.

\section{Nonvanishing of periodic families}\label{sec:detection}
The previous section gave us all of the tools we need to prove the rest of our nonvanishing statements. These are a consequence of the idea that the $E_2$-page of the ANSS for $\Sph$ in filtration $\leq 2$ is quite explicit from work of Novikov \cite{anoneline} and Miller--Ravenel--Wilson \cite{MRW77}. We begin in filtration $\leq 1$.

\begin{lemma}\label{lm:sourcehaslowfiltration}
    Let $\varphi\colon \1 \to X$ be a map of synthetic spectra and $x\in \pi_{s,f}\1$. If $\varphi(x)$ is nonzero and $\tau^{r}$-torsion but not $\tau^{r-1}$-torsion for some $r\geq f-2$, then $x$ is $\tau$-torsion free.
\end{lemma}

\begin{proof}
    By assumption, $x$ is either $\tau^{s}$-torsion for some $s\geq r$ or $\tau$-torsion free. As $r\geq f-2$, we see that $x$ cannot be $\tau^{s}$-torsion for any $s\geq r$, as then it would be hit by a $d_{s+1}$-differential with source in filtration $\leq 1$, and all such classes are permanent cycles by \cite{anoneline}; also see \cite[Th.2.2]{MRW77}. In particular, we see $x$ is $\tau$-torsion free.
\end{proof}

This immediately gives half of our nonvanishing statements.

\begin{prop}\label{pr:nonvanishingbeta1beta63}
    Choose elements $x_i \in \{\be_{1+9s}, \be_{6+9s/3}\}_{s\geq 0} \subseteq \pi_{\ast,2} \1$. Then for all $t\geq 0$:
    \begin{enumerate}
        \item $\prod_{i=1}^5 x_i$ is $\tau^8$-torsion in $\TMF_\BP$ and $\J^{2}_\BP$, and $\tau$-torsion free in $\1$.
        \item $[\al_1\be_{3+9t/3}] x_i$ and $[\al_1\be_{7+9t}] x_i$ are $\tau^4$-torsion in $\TMF_\BP$ and $\tau$-torsion free in $\J^{2}_\BP$.
        \item $\be_{2+9t} x_i$ and $\be_{5+9t} x_i$ are $\tau$-torsion free in $\J^{2}_\BP$.
        \item $\be_{2+9t} \prod_{i=1}^2 x_i$ and $\be_{5+9t} \prod_{i=1}^2 x_i$ are $\tau^4$-torsion in $\J^{2}_\BP$ and $\tau$-torsion free in $\1$.
    \end{enumerate}
\end{prop}

\begin{proof}
    For part 1, we appeal to \Cref{lm:sourcehaslowfiltration} and the spectral sequence of \Cref{pr:manssforJ}, and observe that the $d_9$-differential hitting such a product has source of filtration $1$. 
    For part 2, 
    these products are hit by $d_5$-differentials in the ANSS for $\TMF$ which do not lift to differentials in $\J^2_\BP$ by \Cref{pr:manssforJ}. Part 3 and the $\tau^4$-torsion statement of part 4 both follow from \Cref{pr:manssforJ}. As the products in part 4 have filtration 6 and are hit by $d_5$-differentials in the spectral sequence associated with $\J^2_\BP$, see \Cref{pr:manssforJ}, then by \Cref{lm:sourcehaslowfiltration} we see that these products are $\tau$-torsion free in $\1$.
\end{proof}

For the rest of our nonvanishing statements, we want a refinement of \Cref{lm:sourcehaslowfiltration} stating that certain elements of Adams--Novikov filtration 2 are not detected by $\J^2_\BP/\tau$.

\begin{prop}\label{lm:filtrationtwo_stepping_stone}
	Let $n\geq 0$ and $3\nmid s\geq 1$, $i\geq 0$, and $j\geq 1$ be such that $\beta_{3^n s/j,i+1}$ exists in $\pi_{\ast,\ast} \1/\tau$; see \cite[Th.2.6]{MRW77}. If $n\geq 2$, then the image of $\beta_{3^n s/j,i+1}$ in $\pi_{\ast,\ast} \J^2_\BP/\tau$ is zero.
\end{prop}

\begin{cor}\label{lm:filtrationtwo}
	The mod $\tau$-synthetic unit map $\pi_{s,f} \1/\tau \to \pi_{s,f} \J^2_\BP/\tau$ is zero for $f=2$ and $s\equiv 2, 50, 98, 130$ modulo $144$.
\end{cor}

The cases where $s\equiv 50$ and $130$ are of most interest in this article; the other two cases are included even though they appear in \cite[Pr.3.34]{heighttwojat3}, as we find our proof here simpler.

\begin{proof}
	This whole proof relies on \cite[Th.2.6]{MRW77}. The only classes in the ANSS for $\Sph$ in filtration 2 are divided $\beta$-family classes of the form $\beta_{3^n s/ j, i+1}$ of degree $2^4 3^n s - 4j - 2$. As $\beta_{3^n s/ j, i+2} = 3\beta_{3^n s/ j, i+1}$ when these classes exist, we set $i=0$ without loss of generality. If $n=0$, then $j = 1$, hence the degree of $\beta_{s/ 1}$ is $16s-6$. This degree is congruent to $10$ modulo $16$, and all of the candidate degrees are congruent to $2$ modulo $16$, so $n\neq 0$. If $n=1$, then $1\leq j \leq 3$, and $\beta_{3s/j}$ has degree $48s-4j-2$ which is congruent to $34$, $38$, or $42$ if $s\equiv 1$ modulo $3$, and $82$, $86$, or $90$ if $s\equiv 2$ modulo $3$; these congruence classes differ from all of our candidate degrees. For $n\geq2$, we now appeal to \Cref{lm:filtrationtwo_stepping_stone}, which finishes the proof.
\end{proof}

As promised above, we can compute these statements with \Cref{pr:modtau_hurewicz_image} to demonstrate a tight relationship between divided $\beta$-family elements and $\J^2_\BP/\tau$.

\begin{cor}\label{cor:classification}
    Let $n,s,j,i$ be positive integers such that the class $\beta_{3^ns/j,i+1}$ exists in $\pi_{\ast,2} \1/\tau$; see \cite[Th.2.6]{MRW77}. Then $\beta_{3^ns/j,i}$ is nonzero in $\J^2_\BP/\tau$ if and only if $i=0$ and either $n=0$ and $j=1$, or $n=1$ and $j=3$.
\end{cor}

\begin{proof}
    By \Cref{lm:filtrationtwo_stepping_stone}, we see $\beta_{3^ns/j,i}$ is sent to zero for $n\geq 2$. If $n=0$, we have detection by \Cref{pr:modtau_hurewicz_image} and \cite[Th.2.6]{MRW77} forces $j=1$. If $n=1$, we have detection by \Cref{pr:synthetihurewicz,pr:modtau_hurewicz_image} for $j=3$ and for $j=1,2$ these classes are sent to zero for degree reasons. By observation, whenever a divided $\beta$-family element is detected in $\J^2_\BP/\tau$, we see that $i=0$ as all of these classes are 3-torsion.
\end{proof}

To illustrate the idea behind \Cref{lm:filtrationtwo_stepping_stone} before its proof, let us review the non-synthetic situation for $\TMF$, which is also a formal corollary of \Cref{lm:filtrationtwo_stepping_stone}. The class $\Delta^6 \in \pi_{144} \TMF$ reduces to $v_2^9$ modulo $(3,v_1)$ up to a unit. This means that $q_1(q_0(\Delta^{6N}))$ is equal to $v_2^{9N}$ up to a unit in $\pi_{144N} \TMF/(3,v_1)$. As a consequence, $\TMF$ cannot detect the divided $\beta$-family classes $\beta_{i}$ if $i$ is divisible by 9, as one explicitly has the computation
\[\beta_{i} = \partial_0(\partial_1(v_2^{i})) \ueq \partial_0( \partial_1(q_1(q_0(\Delta^{2i/3}))) ) = 0 \in \pi_{16i-6} \TMF,\]
where $\ueq$ denotes equality up to a unit. This argument is moot, as the group $\pi_{16i-6} \TMF$ vanishes in this case. Nevertheless, it is a template for the proof of the more general statement \Cref{lm:filtrationtwo_stepping_stone} for $\J_\BP^2$. The first step is a lift of $v_2^9$ in $\J^2_\BP/27$.

\begin{lemma}\label{nicer_class}
	There is a unit $V \in \pi_{144,0} \J^2_\BP/27$ such that the image of $V$ in $\J^2_\BP/(\tau, 3,v_1^{3})$ is equal to $v_2^9$ up to a unit, and the image of $V$ in $\TMF_\BP/27$ is equal to the reduction of $\Delta^6$. In particular, as $\J^2_\BP/27$ has an $\E_2$-structure by (the proof of) \cite[Th.1.2]{burklundmultiplicativestructuresmoorespectra}, it follows that all $\J^2_\BP/27$-modules, such as $\TMF_\BP/(\tau,3)$ and $\J^2_\BP/(\tau,3)$, are $V$-local.
\end{lemma}

To make sense of the first point above, recall that the right unit formula from (\ref{eq:right_unit}) shows that $v_2^3$ is primitive modulo $(3,v_1^3)$, giving a chosen class $v_2^3 \in \pi_{48,0} \1/(\tau,3,v_1^3)$. A class playing the same role as $V$ was used in \cite[Pr.3.17]{heighttwojat3} and denoted as $\widetilde{\partial \Delta^6}$, but as mentioned above, we find our approach here simpler. 

\begin{proof}
	It is easy to check that the horizontal maps in the commutative diagram
	\[\begin{tikzcd}
	{\pi_{144,0} \J^2_\BP/27}\ar[d]\ar[r, "{\simeq}"]	&	{\pi_{144,0} \TMF/27 \simeq \Z/27\{c_4^{3r}\Delta^{6-r}\}_{r\geq0}}\ar[d]	\\
	{\pi_{144,0} \J^2_\BP/(3,v_1^{3})}\ar[r, "{\simeq}"]	&	{\pi_{144,0} \TMF/(3,v_1^{3}) \simeq \F_3\{\Delta^{6}\}}
	\end{tikzcd}\]
	are isomorphisms using \Cref{pr:manssforJ}, and that the vertical maps are surjective. Moreover, the mod $\tau$-reduction maps from all of these groups are isomorphisms. This reduces us to finding a class in $\pi_{144,0} \TMF_\BP/27$ congruent to $v_2^9$ modulo $(3,v_1^3)$. By \cite[Lm.6.2]{tmfthree}, the class $v_2^9$ is congruent to $-\Delta^6$ modulo $(3,v_1^3)$, so $\Delta^6$ is our desired class. The fact that $V$ is unit follows from the fact that it maps to the unit $\Delta^6$ in $\TMF/27$ and the 5-lemma from the defining fibre sequence for $\J^2_\BP$, or alternatively one can produce its inverse $V^{-1}$ by lifting the reduction of $\Delta^{-6}$ in $\pi_{-144,0} \TMF_\BP/27$.
\end{proof}

\begin{proof}[Proof of \Cref{lm:filtrationtwo_stepping_stone}]
	Recall that $\beta_{3^n s/ j,i+1}$ can be written as $\partial_{0^{i+1}}(\partial_{1^j}(v_2^{3^n s}))$ modulo $\tau$, for various $n,s,j,i$ as prescribed by \cite[Th.2.6]{MRW77}. To show that these classes vanish, it suffices, as in the proof of \Cref{lm:filtrationtwo}, to take $i=0$. In fact, we will show that the image of $\partial_{1^j}(v_2^{3^n s})$ is zero in $\pi_{\ast,1} \J^2_\BP/(\tau,3)$. Recall that $n\geq 2$.
    
    As $v_2^9$ is congruent to the reduction of the class $V$ from \Cref{nicer_class} modulo $(3,v_1^3)$, then $v_2^{3^ns}$ is congruent to the reduction of $V^{3^{n-2}s}$ modulo $(3,v_1^{3^{n-1}})$; here we use that $\1/(\tau,3,v_1^j)$ has a ring structure as we are working modulo $\tau$. In particular, as long as $j\leq 3^{n-1}$, then $v_2^{3^n s}$ lies in the image of the reduction map induced by $q_{1^j}\colon \J^2_\BP/(\tau,3) \to \J^2_\BP/(\tau,3,v_1^j)$, hence $\partial_{1^j}(v_2^{3^n s})$ vanishes. This reduces us to $j> 3^{n-1}$.
	
    Consider the diagram of abelian groups for some large natural number $N$
	\[\begin{tikzcd}
	{\pi_{2^4 3^n s, 0} \J^2_\BP/(\tau, 3, v_1^j)}\ar[d, "{V^N}"]\ar[r, "{\partial_{1^j}}"]	&		{\pi_{2^4 3^n s - 4j -1, 1} \J^2_\BP/(\tau, 3)}\ar[d, "{V^N}"]
    \\
	{\pi_{2^4 3^n s + 144N, 0} \J^2_\BP/(\tau, 3, v_1^j)}\ar[r, "{\partial_{1^j}}"]	&		{\pi_{2^4 3^n s - 4j +144N -1, 1} \J^2_\BP/(\tau, 3),}
	\end{tikzcd}\]
	which commutes as $\partial_{1^j}$ is $\J^2_\BP/27$-linear. Our goal is to show that $\partial_{1^j}(v_2^{3^ns})$ is zero in the upper-right homotopy group above. From \Cref{nicer_class}, the maps $V^N$ are isomorphisms, so it suffices to show that $\partial_{1^j}(V^N v_2^{3^ns})$ is zero in the lower-right homotopy group. Taking $N$ to be large enough, $V^N$ is congruent to $v_2^{9N}$ modulo $(3,v_1^j)$, meaning that 
    \[\partial_{1^j}(V^N v_2^{3^ns}) = \partial_{1^j}(v_2^{3^ns + 9N}) = 0,\]
    the vanishing a consequence of the paragraph above, as we can again take $N$ to be large enough such that our fixed $j \leq (3^{n} + N)/3$. This finishes the proof.
\end{proof}

Similar to \Cref{pr:nonvanishingbeta1beta63}, we now have our final nonvanishing statement.

\begin{prop}\label{pr:secondlayerusingsecondfiltration}
    The classes of (\ref{eq:line2}) and (\ref{eq:line2.5}) are $\tau^4$-torsion in $\J^2_\BP$, the classes of (\ref{eq:line3}) are $\tau^8$-torsion in $\J^2_\BP$, and all are $\tau$-torsion free in $\1$.
\end{prop}

\begin{proof}
    For the classes in (\ref{eq:line2}), we see by \Cref{pr:synthetihurewicz} that they are nonzero in $\J^2_\BP$ and then use \Cref{pr:manssforJ} to see that they are all $\tau^4$-torsion; the classes $\al_1\be_1^2\be_2$, $\al_1\be_2\be_{6/3}^2$, and their associated periodic families are $\tau^4$-divisible by something hit by a $d_9$-differential, hence these classes are $\tau^4$-torsion. In particular, all of these products are either hit by a $d_5$-differential in the ANSS for $\Sph$ or not. If they are killed by a $d_5$, then this would have to be the same $d_5$ as in the modified ANSS of $\J^2_\BP$, which is impossible as the bidegree for the source of this differential is mapped to zero by \Cref{lm:filtrationtwo}. Hence the classes of (\ref{eq:line2}) are $\tau$-torsion free in $\1$. The same goes for the classes of (\ref{eq:line2.5}). The classes in (\ref{eq:line3}) are $\tau^8$-torsion in $\J^2_\BP$ again by inspecting \Cref{pr:manssforJ}, and the same argument as above shows that they are $\tau$-torsion free in $\1$.
\end{proof}

\begin{remark}\label{rmk:morefamilies?}
    It is not clear to us that a mixture of the first two families of (\ref{eq:line2}) is nonzero in $\pi_\ast \Sph$, ie, the family of the form $\al_1 \beta_{1+9s}\beta_{2+9t}\beta_{6+9w/3}$. Indeed, if we tried to copy the above argument in this case, we would need to use that the classes in $\pi_{122,2} \J_\BP^2/\tau$ do not lie in its mod $\tau$-Hurewicz image, but we know that these classes do lie in the Hurewicz image by \Cref{pr:modtau_hurewicz_image}. On the other hand, mixtures of families from (\ref{eq:line2.5}) such as $\prod_{i=1}^2\beta_{1+9s_i} [\al_1 \beta_{3+9t/3}]$ definitely do vanish, at least for $t=0$, by \Cref{cor:vanishinghalf}.
\end{remark}

\begin{proof}[Proof of \Cref{main:nonvanishing}]
    Parts 1 and 4 of \Cref{pr:nonvanishingbeta1beta63} give us (\ref{eq:line1}), and \Cref{pr:secondlayerusingsecondfiltration} gives us the rest.
\end{proof}

\begin{proof}[Proof of \Cref{main:equivalences}]
    Combine \Cref{cor:vanishinghalf} with \Cref{main:nonvanishing}.
\end{proof}


\addcontentsline{toc}{section}{References}
\scriptsize
\bibliography{references}
\bibliographystyle{alpha}

\end{document}